\documentclass{amsart}

\usepackage{graphicx,amsmath,amssymb,amsthm,geometry,caption}
\usepackage{indentfirst}
\usepackage{amsrefs}

\allowdisplaybreaks

\theoremstyle{plain}
\newtheorem{theorem}{Theorem}[section]

\newtheorem{corollary}[theorem]{Corollary}

\theoremstyle{definition}
\newtheorem*{definition}{Definition}
\newtheorem{example}[theorem]{Example}

\theoremstyle{remark}
\newtheorem{remark}[theorem]{Remark}

\def\*#1{\mathbf{#1}}

\begin{document}

\title[A Note on the Ces\`{a}ro Condition]
{A Note on the Ces\`{a}ro Condition for Curves in the Three-Dimensional Heisenberg Group}

\author[Y.C. Huang]{Yen-Chang Huang}
\address{Yen-Chang Huang \\ Department of Applied Mathematics\\National Yang Ming Chiao Tung University\\Hsinchu City 30010, Taiwan}
\email{ychuang0802@nycu.edu.tw}
\thanks{This work was funded by the National Science and Technology Council (NSTC), Taiwan, with grant number 115-2115-M-A49-001-MY2.}

\subjclass[2020]{Primary: 53A04, 32V30; Secondary: 53C23, 53C17}
\keywords{Ces\`{a}ro condition, adapted Bertrand curves, Heisenberg groups}

\begin{abstract}
We study three consequences of the adapted Frenet--Serret formulas for horizontally regular curves in the three-dimensional Heisenberg group $\mathbb{H}_1$, viewed as the flat model in pseudo-Hermitian geometry. First, we derive the Ces\`{a}ro immobility system and solve it explicitly under the fixed standard coordinate identification
$\mathbb H_1\simeq\mathbb R^3$. The solution yields identities relating the $p$-curvature and contact normality to the radial and vertical coordinates of a curve, and gives a criterion for a curve to lie on a rotationally symmetric surface. We also describe the constant $p$-curvature case and illustrate the formulas on the Euclidean sphere and the Pansu sphere. Second, we classify nontrivial adapted
Bertrand mates through their horizontal projections, derive the
relation between their $p$-curvatures, and describe the remaining
freedom in the vertical component. Finally, we classify curves whose
position vectors lie in the planes spanned by pairs of vectors in the
adapted frame. The three cases yield curves contained in the
$xy$-plane, a vertical plane through the $z$-axis, or a circular
cylinder. These results isolate the effects of the fixed Reeb
direction on several classical constructions for curves in Euclidean
space.
\end{abstract}

\maketitle

\section{Introduction}
Let $\*r$ be a regular curve in the three-dimensional Euclidean space $\mathbb{R}^3$, parametrized by arc length. Its Frenet--Serret frame $\{\*t,\*n,\*b\}$ satisfies
\begin{align*}
\begin{array}{rlll}
\*t'=& &\kappa \*n,& \\
\*n'=&-\kappa \*t & &+\tau \*b, \\
\*b'=&-\tau \*n,&
\end{array}
\end{align*}
where $\kappa$ and $\tau$ are the curvature and torsion. These formulas give several classical characterizations of curves in $\mathbb{R}^3$. For example, the position vector can be used to characterize curves whose position vectors lie in the osculating, normal, or rectifying planes; see \cite[Propositions 4.7, 4.8]{millman1977elements} and Chen \cite{chen2003does}.

In this paper, we study applications of the adapted Frenet--Serret formulas to horizontally regular curves in the three-dimensional Heisenberg group $\mathbb{H}_1$, regarded as a pseudo-Hermitian manifold with vanishing Tanaka--Webster curvature. For basic background on the Heisenberg group from this perspective, we refer to \cites{dong2021schwarz,chang2021sharp} and the references therein. In addition, Chiu--Huang--Lai \cite{CHL} derived an analogous Frenet--Serret formula for the adapted frame of such curves (see \eqref{Frenet} below) and proved the Fundamental Theorem of Curves in $\mathbb{H}_1$. Related aspects of curve geometry in contact and pseudo-Hermitian $3$-manifolds have been studied for slant and Legendre curves by Cho--Lee \cite{cho2008slant}, Lee \cite{lee2010legendre}, and \"Ozg\"ur--G\"uven\c{c} \cite{ozgur2012slant}. Global properties of curves in $\mathbb{H}_1$ and the CR sphere were studied by Chiu--Ho \cite{chiu2019global}. Two features of the adapted frame are particularly important for the questions considered here: the vector $\*b=\*T$ is fixed, and the derivative of the adapted normal vector field contains the term $-\*b$; see the third equation in \eqref{Frenet}. These features make several classical curve-theoretic constructions in $\mathbb{H}_1$ different from their Euclidean counterparts. The purpose of this paper is to describe three such consequences.

First, we derive the Ces\`{a}ro immobility condition \eqref{cesarocondition}. In classical Euclidean geometry, Cesàro's immobility condition describes the variation of the coordinates of a fixed point with respect to a moving Frenet frame and can be used, in particular, to characterize spherical curves in terms of their curvature and torsion; see Cesàro \cite{cesaro1896lezioni}, Lagally \cite{Lagally1927}, Struik \cite{Struik}, and Remark \ref{importantcondition} below. In $\mathbb{H}_1$, the first two equations of the corresponding system, \eqref{cesarocondition}, are independent of the third component; this decoupling allows us to solve explicitly for the first two adapted-frame components and then recover the third. We consequently obtain identities relating the $p$-curvature $\kappa$ and contact normality $\tau$ to the radial and vertical coordinates of a curve. These identities also give a criterion
for a curve to lie on a rotationally symmetric surface; see Theorem \ref{rotation1}. The formulation allows the parameter of the generating profile curve to vary independently of the horizontal arc-length parameter of the curve,
and therefore does not require the profile parameter to coincide with
the horizontal arc-length parameter.

Second, we study Bertrand mates. In $\mathbb{R}^3$, a Bertrand curve admits a mate whose principal normal lines coincide at corresponding points, and the classical theory gives a linear relation between curvature and torsion of the curve; see, for example, \cite[p.~27]{do2016differential}. Bertrand-type curves have also been studied in non-Euclidean ambient spaces and in the theory of framed curves; see \cites{lucas2012bertrand,balgetir2004null,jin2008null,
honda2020bertrand}. Related Bertrand and Mannheim curves have also been considered in Heisenberg-type metric settings: K\"orp\i nar--Turhan \cite{korpinar2011bertrand} studied Bertrand curves in $\mathrm{Heis}^3$, Ceylan--Ergin \cite{ceylan2015bertrand} studied Bertrand mates of biharmonic curves in Cartan--Vranceanu $3$-spaces, and Tul--Sar{\i}o\u{g}lugil \cite{tul2011mannheim} studied Mannheim partner curves in the three-dimensional Heisenberg group. In the setting of three-dimensional Lie groups with a bi-invariant metric, Okuyucu--G\"ok--Yayl{\i}--Ekmekci \cite{okuyucu2017bertrand}
characterized Bertrand curves in terms of the curvature and harmonic curvature. Motivated by these works, we adopt the pseudo-Hermitian viewpoint and use the adapted frame $\{\*t,\*n,\*T\}$ in \eqref{Frenet} to give a necessary and sufficient condition for a given horizontally regular curve $\*r$ in $\mathbb{H}_1$ to admit an adapted Bertrand mate $\bar{\*r}=(\bar x,\bar y,\bar z)$. We also characterize the $\bar x$- and $\bar y$-components of the mate and derive its $p$-curvature $\bar\kappa$; see Theorem \ref{existencethm}. We observe, however, that the $\bar z$-component is not determined by the equality of the adapted normal fields, $\bar{\*n}(\bar s(s))=\*n(s)$. This yields a substantial difference from the classical Euclidean Bertrand theory.

 
Finally, we classify curves according to the position of their position vectors relative to the planes spanned by pairs of vectors in the adapted frame $\{\*t,\*n,\*b\}$. This is a counterpart of the classical position-vector characterizations in $\mathbb{R}^3$; see \cite[Propositions 4.7, 4.8]{millman1977elements} for the osculating and normal cases and Chen \cite{chen2003does} for the rectifying case. If the position vector of a horizontally regular curve lies in $\operatorname{span}\{\*t,\*n\}$, then the curve lies in the $xy$-plane. If it lies in $\operatorname{span}\{\*t,\*b\}$, then the curve lies in a vertical plane through the $z$-axis. If it lies in $\operatorname{span}\{\*n,\*b\}$, then the projection onto the $xy$-plane lies on a circle and the curve lies on a circular cylinder. Its height is affine precisely when its contact normality is constant. In the constant-$\tau$ case, constant height gives a planar circle, whereas nonconstant height gives a nondegenerate circular helix.

\section{Preliminaries}\label{Preliminary}
Chiu--Huang--Lai \cite{CHL} studied the group $PSH(1)$ of pseudo-Hermitian transformations of the Heisenberg group $\mathbb{H}_1$ and obtained
\[
PSH(1)=\mathcal{U}(1)\ltimes\mathcal{T}(1),
\]
the semidirect product of the unitary group $\mathcal{U}(1)$ and the group $\mathcal{T}(1)$ of left translations in $\mathbb{H}_1$. Thus, PSH(1) is precisely the group of rigid motions of $\mathbb{H}_1$ when viewed as a pseudo-Hermitian manifold with vanishing Tanaka--Webster curvature. For $\*p=(x,y,z)$ and $\*q=(a,b,c)$ in $\mathbb{H}_1$, the left translation $L_{\*p}$ is
\[
L_{\*p}(\*q)=(a+x,b+y,c+z+ya-xb).
\]
The standard left-invariant frame at $\*p$ is
\begin{equation}\label{standardbasis}
\begin{aligned}
\mathring{\*e}_1(\*p)&=\frac{\partial}{\partial x}+y\frac{\partial}{\partial z}:=(1,0,y),\\
\mathring{\*e}_2(\*p)&=\frac{\partial}{\partial y}-x\frac{\partial}{\partial z}:=(0,1,-x),\\
\*T(\*p)&=\frac{\partial}{\partial z}:=(0,0,1).
\end{aligned}
\end{equation}
The almost complex structure $J$ on the contact bundle
$\xi=\operatorname{span}\{\mathring{\*e}_1,\mathring{\*e}_2\}$ is
defined by
$J\mathring{\*e}_1=\mathring{\*e}_2$, $J\mathring{\*e}_2=-\mathring{\*e}_1$, and is extended to $T\mathbb{H}_1$ by setting $J\*T=\*0$.

Throughout the paper, all curves are assumed to be smooth and, unless otherwise stated, are parametrized by horizontal arc
length.

On the contact bundle \(\xi\), we use the Levi metric
$h(\*X,\*Y):=\frac12 d\Theta(\*X,J\*Y)$ for any $\*X,\*Y\in\xi$, where $\Theta=dz+xdy -ydx$ is the canonical contact form, and extend \(h\) to the adapted metric on \(T\mathbb H_1\) by
\[
g_\Theta(\*X,\*Y)
:=
\frac12 d\Theta(\*X,J\*Y)+\Theta(\*X)\Theta(\*Y),
\qquad \*X,\*Y\in T\mathbb H_1.
\] We write
\(\langle \*X,\*Y\rangle:=g_\Theta(\*X,\*Y)\). Thus
\(\{\mathring{\*e}_1,\mathring{\*e}_2,\*T\}\) is an orthonormal frame with respect to the metric $\langle \cdot , \cdot \rangle$.

For any curve \(\mathbf r\), its velocity naturally decomposes as
\[
\mathbf r'=\mathbf r'_\xi+\Theta(\mathbf r')\*T,
\]
where \(\mathbf r'_\xi\in\xi\). The curve is called \textit{horizontally regular} if \(\mathbf r'_\xi\) never vanishes along $\*r$. It can then be
reparametrized by horizontal arc length \(s\), so that
$\lvert\mathbf r'_\xi\rvert :=\sqrt{\langle \*r'_\xi, \*r'_\xi \rangle}=1$.

Setting $\mathbf t=\mathbf r'_\xi$ and $\mathbf n=J\mathbf t$, we define
\begin{align}\label{tau}
\kappa(s):=\langle \mathbf t'(s),\mathbf n(s)\rangle,
\qquad
\tau(s):=\langle\mathbf r'(s),\*T(s)\rangle
        =\Theta(\mathbf r'(s)).
\end{align}
The functions $\kappa$ and $\tau$ are called the $p$-curvature and the contact normality, respectively. In \cite[Theorem 1.2]{CHL}, the authors proved the Fundamental Theorem of Curves in $\mathbb{H}_1$, which states that these two invariants determine a horizontally regular curve uniquely up to a pseudo-Hermitian transformation. Higher-dimensional analogues can be found in \cite{chiu2018differential}.
If $\*r(u)=(x(u),y(u),z(u))$ is parametrized by an arbitrary parameter $u$, then
\begin{equation}\label{pcurve}
\begin{aligned}
\kappa(u)&=\frac{x'y''-x''y'}{\big((x')^2+(y')^2\big)^{3/2}},\\
\tau(u)&=\frac{xy'-x'y+z'}{\big((x')^2+(y')^2\big)^{1/2}}.
\end{aligned}
\end{equation}
In particular, $\tau\equiv0$ if and only if $\*r$ is horizontal.

\begin{remark}
Let $\nabla$ denote the Tanaka--Webster connection of the standard pseudo-Hermitian structure on $\mathbb H_1$, and let $D$ denote the standard
flat affine connection under the coordinate identification $\mathbb H_1\simeq\mathbb R^3$. For a horizontally regular curve
$\*r(s)=(x(s),y(s),z(s))$, parametrized by horizontal arc length, 
write
\[
\*t:=\*r'_\xi = x'\mathring{\*e}_1+y'\mathring{\*e}_2,
\qquad
\*n:= J\*r'_\xi =-y'\mathring{\*e}_1+x'\mathring{\*e}_2.
\] Since
$
D_s\mathring{\*e}_1=y'\*T$, $D_s\mathring{\*e}_2=-x'\*T$, we obtain
\[
\begin{aligned}
D_s\*t
&=
x''\mathring{\*e}_1+y''\mathring{\*e}_2
+\bigl(x'y'-y'x'\bigr)\*T \\
&=
x''\mathring{\*e}_1+y''\mathring{\*e}_2.
\end{aligned}
\]
On the other hand, the standard left-invariant frame
$\{\mathring{\*e}_1,\mathring{\*e}_2,\*T\}$ is parallel with respect to
$\nabla$, and hence 
$$\nabla_s\*t
=
x''\mathring{\*e}_1+y''\mathring{\*e}_2. $$ Therefore
$D_s\*t=\nabla_s\*t$,
and consequently
$
\kappa
=
\langle D_s\*t,\*n\rangle
=
\langle\nabla_s\*t,\*n\rangle$. Thus the use of the auxiliary flat connection $D$ gives the same $p$-curvature as the intrinsic Tanaka--Webster  connection. The contact normality $\tau=\Theta(\mathbf r')$ is independent of either
connection.

Notice, however, that $D$ and $\nabla$ do not coincide in general. Indeed,
using $(x')^2+(y')^2=1$,
\[
D_s\*n
=
\nabla_s\*n-\*T.
\]
This accounts for the $-\*T$ term in the adapted Frenet--Serret formula. By \cite[Theorem 3.1]{CHL}, every element of \(PSH(1)\) acts affinely in these coordinates; hence
$\Phi^*D=D$ for every $\Phi\in PSH(1)$. Thus \(D\) is \(PSH(1)\)-invariant.
\end{remark}

For a horizontally regular curve parametrized by horizontal arc length, the adapted Frenet--Serret formulas are
\begin{align}\label{Frenet}
\left\{
\begin{array}{rcl}
\*r'&=&\*t+\tau\*b,\\
\*t'&=&\kappa\*n,\\
\*n'&=&-\kappa\*t-\*b,\\
\*b'&=&\*0,
\end{array}
\right.
\end{align}
where
\begin{equation}\label{basissetting}
\*t=\*r'_\xi,
\qquad
\*n=J\*r'_\xi,
\qquad
\*b=\*T;
\end{equation} see \cite{CHL}. 
If $\*r=(x,y,z)$, then
$\*r'=x'\mathring{\*e}_1+y'\mathring{\*e}_2+(-x'y+xy'+z')\*T$, 
and hence
\begin{align}
\*t&=x'\mathring{\*e}_1+y'\mathring{\*e}_2=(x',y',x'y-xy'),\label{tangent1}\\
\*n&=J\*t=-y'\mathring{\*e}_1+x'\mathring{\*e}_2=(-y',x',-yy'-xx').\label{normal1}
\end{align}

\begin{remark}\label{distance1}
Using the standard coordinate identification
$\mathbb{H}_1\simeq\mathbb{R}^3$, and since 
$(x,y,z)=x\mathring{\*e}_1+y\mathring{\*e}_2+z\*T$, we can express the position vector of $\*r$ in the adapted frame. Indeed, write
\[
\*t=\cos\alpha\,\mathring{\*e}_1+\sin\alpha\,\mathring{\*e}_2,
\qquad
\*n=-\sin\alpha\,\mathring{\*e}_1+\cos\alpha\,\mathring{\*e}_2.
\]
Then
\begin{equation}\label{position}
\*r=\widetilde{x}\,\*t+\widetilde{y}\,\*n+z\,\*b,
\end{equation}
where $\widetilde{x}=x\cos\alpha+y\sin\alpha$ and $\widetilde{y}=-x\sin\alpha+y\cos\alpha$. Consequently,
\begin{equation}\label{distance}
x^2+y^2=\widetilde{x}^{\,2}+\widetilde{y}^{\,2}.
\end{equation}
Thus, the Euclidean distance of $\*r$ to the $z$-axis is $\sqrt{\widetilde{x}^{\,2}+\widetilde{y}^{\,2}}$, and the coefficient of $\*b$ in \eqref{position} is exactly the Euclidean height $z$.
\end{remark}

\section{Ces\`{a}ro conditions}\label{sectioncesaro}
Throughout this section, expressions involving sums of points and adapted-frame vectors are understood through the fixed standard coordinate identification
$
\mathbb{H}_1\simeq\mathbb{R}^3$ 
introduced in Remark \ref{distance1}.

Let $\*r:I\to\mathbb{H}_1$ be a horizontally regular curve parametrized by horizontal arc length. Consider
\begin{align}\label{curve1}
\bar{\*r}(s)=\*r(s)+u_1(s)\*t(s)+u_2(s)\*n(s)+u_3(s)\*b(s).
\end{align}
Differentiating and using \eqref{Frenet}, we obtain
\begin{align}\label{maineq1}
\bar{\*r}'
&=(1+u_1'-\kappa u_2)\*t+(u_2'+\kappa u_1)\*n+(\tau+u_3'-u_2)\*b.
\end{align}
Therefore $\bar{\*r}$ is constant if and only if
\begin{align}\label{cesarocondition}
\left\{
\begin{array}{rl}
u_1'&=\kappa u_2-1,\\
u_2'&=-\kappa u_1,\\
u_3'&=u_2-\tau.
\end{array}
\right.
\end{align}
This is the Ces\`{a}ro immobility condition in $\mathbb{H}_1$.

\begin{remark}\label{importantcondition}
Writing
$
\*r=\widetilde{u}_1\*t+\widetilde{u}_2\*n+\widetilde{u}_3\*b$,
and taking $\bar{\*r}=\*0$ in \eqref{curve1}, we see that $u_i=-\widetilde{u}_i$ satisfy \eqref{cesarocondition}.  In Euclidean space the corresponding Ces\`{a}ro immobility condition is, up to the convention for the sign of the torsion, 
\begin{align}
\left\{
\begin{array}{rl}
u_1'&=\kappa u_2-1,\\
u_2'&=-\kappa u_1+\tau u_3,\\
u_3'&=-\tau u_2.
\end{array}
\right.
\end{align} See Kowalewski \cite[Sec.~3, Eq.~(7), p.~77]{Kowalewski}; see also Lagally \cite{Lagally1927} for the classical moving-frame formulation. The absence of $u_3$ from the first two equations of \eqref{cesarocondition} is a special simplification in the Heisenberg group setting and allows us to solve explicitly for $u_1$ and $u_2$. 
\end{remark}

We now solve \eqref{cesarocondition} directly. Fix $s_0\in I$ and define
\begin{equation}\label{thetaAB}
\theta(s)=\int_{s_0}^s\kappa(t)\,dt,
\qquad
A(s)=\int_{s_0}^s\cos\theta(t)\,dt,
\qquad
B(s)=\int_{s_0}^s\sin\theta(t)\,dt.
\end{equation}
Set $w=u_1+i u_2$. The first two equations in \eqref{cesarocondition} give
$w'+i\kappa w=-1$, so
$\big(e^{i\theta}w\big)'=-e^{i\theta}$.
Hence all solutions are
\begin{align}
 u_1(s)&=(c_1-A(s))\cos\theta(s)+(c_2-B(s))\sin\theta(s),\label{u1solution}\\
 u_2(s)&=-(c_1-A(s))\sin\theta(s)+(c_2-B(s))\cos\theta(s),\label{u2solution}\\
 u_3(s)&=c_3+\int_{s_0}^s\big(u_2(t)-\tau(t)\big)\,dt,\label{u3solution}
\end{align}
where $c_1,c_2,c_3$ are arbitrary constants. These formulas remain valid even when $\kappa$ vanishes. In particular, if $\kappa\equiv0$, then
\[
u_1=c_1-(s-s_0),
\qquad
u_2=c_2,
\qquad
u_3=c_3+c_2(s-s_0)-\int_{s_0}^s\tau(t)\,dt.
\]
On any interval on which $\kappa\neq0$, elimination of $u_2$ and $u_1$, respectively, also gives
\begin{align}
 u_1''-\frac{\kappa'}{\kappa}u_1'+\kappa^2u_1-\frac{\kappa'}{\kappa}&=0,\label{u1eqn}\\
 u_2''-\frac{\kappa'}{\kappa}u_2'+\kappa^2u_2-\kappa&=0.\label{u2eqn}
\end{align}

We next give a precise formulation for curves lying on a surface of revolution.

\begin{theorem}\label{rotation1}
Let \(J\subset\mathbb R\) be an interval and 
\(f,g:J\to\mathbb R\) be smooth functions satisfying $g>0$ and $f'(v)^2+g'(v)^2>0$ for all $v\in J$. Assume the profile curve $v\mapsto (g(v), 0, f(v))$ is an embedding and define the rotationally symmetric surface $\Sigma$ in $\mathbb{H}_1$ about the \(z\)-axis by
\[
\Sigma
=
\bigl\{
(g(v)\cos t,g(v)\sin t,f(v)):
v\in J,\ 0\leq t<2\pi
\bigr\}.
\]
Let \(\mathbf r:I\to\mathbb H_1\) be a horizontally regular curve
parametrized by horizontal arc length. Choose the constants in \eqref{u1solution}--\eqref{u3solution} so that
\begin{equation}\label{originchoice}
\*r(s)+u_1(s)\*t(s)+u_2(s)\*n(s)+u_3(s)\*b(s)=\*0.
\end{equation}
Then $\*r(I)\subset\Sigma$ if and only if there exists a unique function $v:I\to J$ such that
\begin{equation}\label{edistance}
 u_1(s)^2+u_2(s)^2=g(v(s))^2,
\qquad
 u_3(s)=-f(v(s))
\end{equation}
for all $s\in I$.
\end{theorem}

\begin{proof}
By \eqref{originchoice}, we have
$\*r=-u_1\*t-u_2\*n-u_3\*b$. Remark \ref{distance1} therefore gives
\[
\sqrt{x(s)^2+y(s)^2}=\sqrt{u_1(s)^2+u_2(s)^2},
\qquad
z(s)=-u_3(s).
\]
A point $(x,y,z)$ lies on $\Sigma$ exactly when there is a profile parameter $v\in J$ such that $x^2+y^2=g(v)^2$ and $z=f(v)$. The equivalence follows.
\end{proof}

Notice that the parameter $v$ of the generating profile curve $v\mapsto (g(v), 0, f(v))$ should be regarded as independent of the horizontal arc-length parameter $s$ of the curve. In general, one only requires the existence of a reparametrization $v=v(s)$ such that the relations below in \eqref{rzcondition} hold, rather than imposing $v=s$.

The next identities are useful independently of a prescribed surface.

\begin{corollary}\label{nscondition1}
Let $\*r(s)=(x(s),y(s),z(s))$ be a horizontally regular curve parametrized by horizontal arc length, with $\kappa(s)\neq0$, and define the radial and height functions, respectively, by
\begin{align}\label{radialheight}
R(s)=\sqrt{x(s)^2+y(s)^2},
\qquad
Z(s)=z(s).
\end{align}
Then
\begin{align}
 Z'-\tau&=\frac{\frac12(R^2)''-1}{\kappa},\label{coro1}\\
 Z''-\tau'&=-\frac{\kappa}{2}(R^2)'.\label{coro2}
\end{align}
If, in addition, $\*r\subset\Sigma$ as in Theorem \ref{rotation1} with profile parameter $v(s)$, then the same identities hold with
\begin{align}\label{rzcondition}
R(s)=g(v(s)),
\qquad
Z(s)=f(v(s)).
\end{align}
\end{corollary}

\begin{proof}
Choose $u_i$ as in \eqref{originchoice}. By \eqref{originchoice} and Remark \ref{distance1},
we have $u_1^2+u_2^2=R^2$ and 
$u_3=-Z$. Differentiating the first identity and using the first two equations of \eqref{cesarocondition}, we get
\[
\frac12(R^2)'=u_1u_1'+u_2u_2'=-u_1,
\]
so
\begin{equation}\label{u1radial}
 u_1=-\frac12(R^2)'.
\end{equation}
The third Ces\`{a}ro equation gives $-Z'=u_2-\tau$, and hence $u_2=-Z'+\tau$. Substituting these expressions into $u_1'=\kappa u_2-1$ gives \eqref{coro1}, while differentiating $u_2=-Z'+\tau$ and using $u_2'=-\kappa u_1$ gives \eqref{coro2}.
\end{proof}

\begin{corollary}\label{nscondition2}
Assume that $\*r$ has nonzero constant $p$-curvature
$\kappa=\kappa_0$, with $R(s)$ and $Z(s)$ defined as in \eqref{radialheight}. Then there are constants $C_1,C_2,C_3,C_4$ such that
\begin{align}
R(s)^2
&=\frac{1}{\kappa_0}\big(-C_1\cos(\kappa_0s)+C_2\sin(\kappa_0s)+C_3\big),\label{knonconstzeroR}\\
Z(s)
&=\int_{s_0}^s\tau(t)\,dt
+\frac{C_1}{2\kappa_0}\sin(\kappa_0s)
+\frac{C_2}{2\kappa_0}\cos(\kappa_0s)
-\frac{s}{\kappa_0}+C_4.\label{knonconstzeroZ}
\end{align}
The constants satisfy the compatibility condition
\begin{align}\label{c3comp}
C_3=
\frac1{\kappa_0}
+\frac{\kappa_0}{4}\bigl(C_1^2+C_2^2\bigr),
\end{align}
and there is no further algebraic compatibility condition among $C_1$, $C_2$, and $C_4$. Under this
condition, the right-hand side of \eqref{knonconstzeroR} is automatically
nonnegative.\end{corollary}

\begin{proof}
Since $\kappa_0$ is constant, differentiating \eqref{coro1} and comparing with \eqref{coro2} yields
\[
(R^2)'''=-\kappa_0^2(R^2)'.
\]
Thus, if $G=(R^2)'$, then $G''+\kappa_0^2G=0$, so
\[
G=C_1\sin(\kappa_0s)+C_2\cos(\kappa_0s).
\]
Integration gives \eqref{knonconstzeroR}, after renaming the integration constant. Substitution into \eqref{coro1} and one further integration gives \eqref{knonconstzeroZ}; the additive constant in $Z$ is independent of the additive constant in $R^2$.

It remains to determine the compatibility condition among the integration constants. By \eqref{u1radial},
\[
u_1
=
-\frac12(R^2)'
=
-\frac12\left(
C_1\sin(\kappa_0s)+C_2\cos(\kappa_0s)
\right).
\]
Since
$u_1'=\kappa_0u_2-1$,
we obtain
$$u_2
=
\frac1{\kappa_0}
-\frac12C_1\cos(\kappa_0s)
+\frac12C_2\sin(\kappa_0s).$$ Hence
\begin{align*}
R^2
&=u_1^2+u_2^2\\
&=
\frac1{\kappa_0^2}
+\frac14(C_1^2+C_2^2)
-\frac{C_1}{\kappa_0}\cos(\kappa_0s)
+\frac{C_2}{\kappa_0}\sin(\kappa_0s).
\end{align*}
Comparing this expression with \eqref{knonconstzeroR} yields \eqref{c3comp}.

To verify the final assertion, let \(C_1,C_2\in\mathbb R\) be
arbitrary. Substituting the compatibility relation into the
right-hand side of \eqref{knonconstzeroR}, we obtain
\[
\begin{aligned}
R^2
&=
\frac{1}{\kappa_0^2}
+\frac{C_1^2+C_2^2}{4}
+\frac{-C_1\cos(\kappa_0s)+C_2\sin(\kappa_0s)}
       {\kappa_0} \\
&=
\frac14
\bigl(C_1\sin(\kappa_0s)+C_2\cos(\kappa_0s)\bigr)^2 \\
&\quad+
\left(
\frac1{\kappa_0}
-\frac{C_1}{2}\cos(\kappa_0s)
+\frac{C_2}{2}\sin(\kappa_0s)
\right)^2
\geq 0.
\end{aligned}
\]
Thus, this relation imposes no additional nonnegativity
restriction on \(C_1\) and \(C_2\). As noted above, \(C_4\) is an
independent additive constant.
\end{proof}

\begin{example}\label{example1}
Let $\mathbb{S}^2_a=\{x^2+y^2+z^2=a^2\}$ be the standard Euclidean sphere, $a>0$. Horizontal curves on $\mathbb{S}^2_a$ do exist. Indeed, in cylindrical coordinates $(x,y)=(R\cos\phi,R\sin\phi)$, the conditions that the curve lie on $\mathbb S_a^2$ and be horizontal are
\begin{equation}\label{spheresystem}
R^2+Z^2=a^2,
\qquad
Z'+R^2\phi'=0.
\end{equation}
If the curve is parametrized by horizontal arc length, then
$
R'^2+R^2\phi'^2=1$. Away from the axis $R=0$, equations \eqref{spheresystem} give, for either choice of orientation,
\[
\phi'=\pm\frac{1}{R\sqrt{1+Z^2}},
\qquad
R'=\pm\frac{Z}{\sqrt{1+Z^2}},
\qquad
Z'=\mp\frac{R}{\sqrt{1+Z^2}}.
\]
Thus there are nontrivial local horizontal curves on the Euclidean sphere. This also illustrates why the profile parameter of a surface of revolution cannot be identified a priori with the horizontal arc-length parameter of an arbitrary curve on the surface.
\end{example}

\begin{example}\label{pansuexample}
For $\lambda>0$, the Pansu sphere $P_\lambda$ is the union of the graphs of $F$ and $-F$, where
\begin{align}\label{pansu1}
F(x,y)=\frac{1}{2\lambda^2}\left(
\lambda r\sqrt{1-\lambda^2r^2}+\arccos(\lambda r)
\right),
\qquad
r=\sqrt{x^2+y^2}\leq\frac1\lambda;
\end{align}
see, for instance, \cite[Lemma~2.4 and Eqs.~(2.2)--(2.4)]{chiu2025invariants}.

A horizontal geodesic joining the two poles is
\begin{align*}
x(s)&=\frac{1}{2\lambda}\sin(2\lambda s),\\
y(s)&=-\frac{1}{2\lambda}\cos(2\lambda s)+\frac{1}{2\lambda},\\
z(s)&=\frac{1}{4\lambda^2}\sin(2\lambda s)-\frac{s}{2\lambda}+\frac{\pi}{4\lambda^2},
\qquad 0\leq s\leq\frac{\pi}{\lambda}.
\end{align*}
A direct computation using \eqref{pcurve} gives
$\kappa=2\lambda$,
$\tau=0$,
and
\[
R(s)=\sqrt{x(s)^2+y(s)^2}=\frac{\sin(\lambda s)}{\lambda}.
\]
Therefore the rotational profile traced by this geodesic is
\begin{equation}\label{pansuprofile}
R(s)=\frac{\sin(\lambda s)}{\lambda},
\qquad
Z(s)=\frac{\sin(2\lambda s)}{4\lambda^2}-\frac{s}{2\lambda}+\frac{\pi}{4\lambda^2}.
\end{equation}
These functions satisfy \eqref{coro1}--\eqref{coro2} with $\kappa=2\lambda$ and $\tau=0$. They also have the form in Corollary \ref{nscondition2}: one may take
\[
C_1=\frac1\lambda,
\qquad
C_2=0,
\qquad
C_3=\frac1\lambda,
\qquad
C_4=\frac{\pi}{4\lambda^2}.
\]
Rotating \eqref{pansuprofile} about the $z$-axis recovers $P_\lambda$.
\end{example}

\section{Bertrand mates}\label{Bertrandmates}
Recall that, in classical Euclidean geometry, a regular curve $\*r(s)$ in $\mathbb{R}^3$ is called a Bertrand curve if it admits a mate $\bar{\*r}(\bar{s})$ whose principal normal line at $\bar{\*r}(\bar s(s))$ coincides with that of $\*r$ at $\*r (s)$ for every $s$. Motivated by this classical notion, we introduce an analogous definition for horizontally regular curves in $\mathbb{H}_1$.

\begin{definition}
A horizontally regular curve $\*r$ in $\mathbb{H}_1$ is called an \emph{adapted Bertrand curve} if there exists a distinct horizontally regular curve $\bar{\*r}$ and an orientation-preserving correspondence $\bar s=\bar s(s)$ such that, under the fixed standard coordinate identification,
\begin{align}\label{betrandcod1}
\bar{\*n}(\bar s(s))=\*n(s)
\quad\text{and}\quad
\pi_{xy}(\bar{\*r}(\bar s(s)))
\neq
\pi_{xy}(\*r(s))
\end{align}
for all $s$. This curve $\bar{\*r}$ is called an \emph{adapted Bertrand mate} of $\*r$.
\end{definition}

\begin{remark}\label{trivialexample}
For a curve parametrized by horizontal arc length, the adapted normal field in \eqref{normalcondi} depends only on $x,y$ and their first derivatives, and is independent of the vertical component $z$. Hence, for every smooth function $h\not\equiv 0$, the horizontally regular curve
$$
\bar{\*r}(s)=\bigl(x(s),y(s),z(s)+h(s)\bigr)
$$
has the same horizontal projection and adapted normal field as $\*r$, while its contact normality satisfies $\bar\tau=\tau+h'$. Thus equality of the adapted normal fields constrains the relation between the horizontal projections but leaves the vertical component undetermined. These trivial vertical modifications correspond to the case $a=0$ in the statement of Theorem~\ref{existencethm} and are excluded by the second condition in \eqref{betrandcod1}, which requires the horizontal projections to be distinct.
\end{remark}

Next, we give a necessary and sufficient condition for a curve to be a nontrivial adapted Bertrand curve.

\begin{theorem}\label{existencethm}
Let $\*r(s)=(x(s),y(s),z(s))$ be a horizontally regular curve parametrized by horizontal arc length. Under the standard coordinate identification $\mathbb{H}_1\simeq\mathbb{R}^3$, a distinct horizontally regular curve
\[
\bar{\*r}(\bar s)
=
(\bar x(\bar s),\bar y(\bar s),\bar z(\bar s))
\]
is an adapted Bertrand mate of $\*r$ if and only if there exists a constant $a\in\mathbb{R}\setminus \{0\}$ such that
\[
1-a\kappa(s)>0
\]
and
\begin{align}\label{xyrelation}
\bar x(\bar s(s))=x(s)-a y'(s),
\qquad
\bar y(\bar s(s))=y(s)+a x'(s),
\end{align}
where
\[
\frac{d\bar s}{ds}=1-a\kappa(s).
\]
Moreover, the $p$-curvature of the adapted Bertrand mate $\bar{\*r}$ is given by 
\[
\bar\kappa(\bar s(s))
=
\frac{\kappa(s)}{1-a\kappa(s)}.
\]
\end{theorem}

\begin{proof}
Suppose first that \eqref{betrandcod1} holds. Since the parameter correspondence is orientation-preserving, set $q(s)=\frac{d\bar s}{ds}>0$. By \eqref{normal1}, we have
\begin{align}\label{normalcondi}
\*n=(-y',x',-yy'-xx')
\text{   and   }
\bar{\*n}
=
\left(
-\bar y_{\bar s},
\bar x_{\bar s},
-\bar y\,\bar y_{\bar s}
-\bar x\,\bar x_{\bar s}
\right).
\end{align}
Here the subscript $\bar s$ denotes differentiation with respect to $\bar s$. The first two coordinate components of
$\bar{\*n}(\bar s(s))=\*n(s)$
give
\begin{align}\label{xycondi}
\bar x_{\bar s}=x' \text{  and  }
\bar y_{\bar s}=y';
\end{align}
the equality of the third coordinate components gives
\begin{align}\label{zcondi}
-\bar y\,\bar y_{\bar s}
-\bar x\,\bar x_{\bar s}
=
-yy'-xx'.
\end{align}
Substituting \eqref{xycondi} into \eqref{zcondi}, we obtain
\[
\bigl(\bar x(\bar s(s))-x(s)\bigr)x'(s)
+
\bigl(\bar y(\bar s(s))-y(s)\bigr)y'(s)
=0.
\]

Define the horizontal displacement vector
\[
\Delta_h(s)
=
\bigl(
\bar x(\bar s(s))-x(s),
\bar y(\bar s(s))-y(s)
\bigr).
\]
Then $\Delta_h$ satisfies 
\begin{align}\label{perpend}
\Delta_h(s)\cdot(x'(s),y'(s))=0.
\end{align}
Since $\*r$ is parametrized by horizontal arc length,
$(x')^2+(y')^2=1$, so $(x',y')$ is a unit vector field. Moreover,
$(x',y')$ and $(-y',x')$ form an orthonormal basis of $\mathbb{R}^2$.
It therefore follows from \eqref{perpend} that there exists a unique
function $a(s)$ such that
\begin{align*}
\Delta_h(s)=a(s)(-y'(s),x'(s)).
\end{align*}
Differentiating this equation with respect to $s$ gives
\[
\Delta_h'
=
a'(-y',x')
+
a(-y'',x'').
\]
Since $(x'',y'')=\kappa(-y',x')$ by the definition of $\kappa$,
substituting this relation into the preceding equation yields
\begin{align}\label{nablacondi1}
\Delta_h'
=
a'(-y',x')
-a\kappa(x',y').
\end{align}

On the other hand, by the definition of $\Delta_h$, together with
\eqref{xycondi} and $q=\frac{d\bar s}{ds}$, we obtain
\begin{align}\label{nablacondi2}
\begin{aligned}
\Delta_h'
&=
\left(
\bar x_{\bar s}\frac{d\bar s}{ds}-x',
\bar y_{\bar s}\frac{d\bar s}{ds}-y'
\right)\\
&=
(q-1)(x',y').
\end{aligned}
\end{align}
Comparing the coefficients in \eqref{nablacondi1} and
\eqref{nablacondi2} with respect to the orthonormal basis
$\{(x',y'),(-y',x')\}$
yields
$a'=0$ and $q=1-a\kappa$. Therefore, $a$ is constant and
\[
\bar x(\bar s(s))=x(s)-a y'(s),
\qquad
\bar y(\bar s(s))=y(s)+a x'(s),
\]
with $\frac{d\bar s}{ds}=1-a\kappa(s)>0$, and the result follows.

Next, we derive the curvature formula. From
\eqref{xycondi} and $q= \frac{d\bar s}{ds}>0$
we obtain
\[
(\bar x_{\bar s\bar s}, \bar y_{\bar s\bar s})
=\frac{1}{q}(x'',y'')
=\frac{\kappa}{q}(-y',x').
\]
Since
$\frac{d}{d\bar s}=\frac{1}{q}\frac{d}{ds}$ and 
$(\bar x_{\bar s},\bar y_{\bar s})=(x',y')$, the definition of $p$-curvature gives
\[
\begin{aligned}
\bar\kappa(\bar s(s))
&=\bar x_{\bar s}\bar y_{\bar s\bar s}
  -\bar y_{\bar s}\bar x_{\bar s\bar s}\\
&=\frac{x'y''-y'x''}{q}
=\frac{\kappa(s)}{q(s)}
=\frac{\kappa(s)}{1-a\kappa(s)}.
\end{aligned}
\]

Conversely, suppose that a constant $a$ satisfies
$q=1-a\kappa>0$, that $d\bar s/ds=q$, and that
\eqref{xyrelation} holds. Since
$x''=-\kappa y'$ and $y''=\kappa x'$, 
differentiating \eqref{xyrelation} gives
\[
\frac{d}{ds}\bar x(\bar s(s))=q x', \qquad
\frac{d}{ds}\bar y(\bar s(s))=q y'.
\]
Consequently, $\bar x_{\bar s}=x'$ and $\bar y_{\bar s}=y'$. Moreover,
\begin{align*}
-\bar y\,\bar y_{\bar s}-\bar x\,\bar x_{\bar s}
&=-(y+ax')y'-(x-ay')x'\\
&=-yy'-xx'.
\end{align*}

All three coordinate components of the adapted normal fields therefore
agree, so $\bar{\*n}(\bar s(s))=\*n(s)$. Hence $\bar{\*r}$ is an adapted Bertrand mate of $\*r$. The preceding computation also proves the asserted curvature formula.

\end{proof}



Next, we show that the Euclidean distance between the horizontal projections of an adapted Bertrand curve and its mate is constant at corresponding points.

\begin{corollary}\label{horizontal-distance}
Let $\bar{\*r}$ be an adapted Bertrand mate corresponding to the constant $a$ in Theorem \ref{existencethm}. Then
\[
\left\|
\pi_{xy}(\bar{\*r}(\bar s(s)))
-\pi_{xy}(\*r(s))
\right\|_{\mathbb R^2}
=
|a|,\]
where $\pi_{xy}$ denotes the projection onto the $xy$-plane.
\end{corollary}

\begin{proof}
By \eqref{xyrelation} in Theorem \ref{existencethm}, we have
\[
\begin{aligned}
\pi_{xy}(\bar{\*r}(\bar s(s)))
-
\pi_{xy}(\*r(s))
&=
\bigl(
\bar x(\bar s(s))-x(s),
\bar y(\bar s(s))-y(s)
\bigr)\\
&=
a(-y'(s),x'(s)).
\end{aligned}
\]
Taking the Euclidean norm in $\mathbb{R}^2$, we obtain
\[
\left\|
\pi_{xy}(\bar{\*r}(\bar s(s)))
-
\pi_{xy}(\*r(s))
\right\|_{\mathbb{R}^2}
=
|a|\sqrt{(y'(s))^2+(x'(s))^2}
=
|a|,
\]
where, in the last equality, we have used the assumption that $\*r$ is parametrized by horizontal arc length.
\end{proof}

Finally, we emphasize that the classical notion of a Mannheim curve does not admit a direct analogue in the geometry of $\mathbb{H}_1$ considered here. Recall that a Mannheim curve in $\mathbb{R}^3$ is a regular curve admitting a mate such that the principal normal line of one curve coincides with the binormal line of the other at corresponding points; see, for instance, \cites{choi2013mannheim,blum1966remarkable,ersoy2012generalized}. In the present adapted setting, an analogous relation involving $\*n$ and $\bar{\*b}$ cannot occur. Indeed, the natural analogue of the binormal vector in the present setting is
$\bar{\*b}=\*T=(0,0,1)$, whereas
$\*n=(-y',x',-yy'-xx')$ has horizontal component $(-y',x')$ of unit length. Hence $\bar{\*b}$ and $\*n$ cannot be parallel. This observation excludes the direct analogue considered above, but it does not rule out all possible Mannheim-type relations in $\mathbb{H}_1$.

\section{A classification of curves in $\mathbb{H}_1$}\label{classification}
We now classify horizontally regular curves according to the position of their position vectors relative to the planes spanned by pairs of vectors in the adapted frame and characterize the corresponding $p$-curvatures.

\begin{theorem}\label{classificationthm}
Let $\*r(s)$ be a horizontally regular curve parametrized by horizontal arc length.
\begin{enumerate}
\item If $\*r(s)\in\operatorname{span}\{\*t(s),\*n(s)\}$ for all $s$, then $\*r$ lies in the $xy$-plane. More precisely:
\begin{enumerate}
\item if $\kappa\equiv0$, then $\*r$ is contained in a straight line in the $xy$-plane;
\item if $\kappa\not\equiv0$, then $\*r$ is a nonlinear planar curve in the $xy$-plane.
\end{enumerate}
\item If $\*r(s)\in\operatorname{span}\{\*t(s),\*b(s)\}$ for all $s$, then $\kappa\equiv0$ and $\*r$ lies in a vertical plane through the $z$-axis.
\item If $\*r(s)\in\operatorname{span}\{\*n(s),\*b(s)\}$ for all $s$, then $\kappa$ is a nonzero constant, the projection of $\*r$ onto the $xy$-plane lies on a circle, and $\*r$ lies on a circular cylinder about the $z$-axis. Moreover, the height $z(s)=A s+B$ for some
constants $A, B \in\mathbb R$ if and only if $\tau$ is constant. In the constant-$\tau$ case, if $z$ is constant, then the curve is a planar circle; if $z$ is nonconstant, then it is a nondegenerate circular helix.
\end{enumerate}
\end{theorem}

\begin{proof}
(1) By \eqref{position} in Remark \ref{distance1}, the coefficient of $\*b$ in the adapted-frame decomposition of the position vector is exactly the Euclidean height $z$ of $\*r$. Therefore, $\*r(s)\in\operatorname{span}\{\*t(s),\*n(s)\}$ for all $s$ implies that the $\*b$-component vanishes, that is, $z(s)=0$ for all $s$. Hence $\*r$ lies in the $xy$-plane.

If $\kappa\equiv0$, its planar projection has zero Euclidean curvature and hence $\*r$ is a straight line. If $\kappa\not\equiv0$, the curve cannot be a straight line, since a straight line has identically vanishing $p$-curvature.

(2) Write $\*r=u_1\*t+u_3\*b$. Differentiating $\*r$ and using \eqref{Frenet} together with the assumption, we obtain
$$
u_1'=1,
\qquad
\kappa u_1=0,
\qquad
u_3'=\tau.
$$
Since $u_1=s+c_1$, the identity $\kappa u_1=0$ implies that $\kappa(s)=0$ whenever $s\neq -c_1$. By continuity of $\kappa$, it follows that $\kappa\equiv0$ on the entire interval. Hence
$$
\*r=(s+c_1)\*t+\left(c_4+\int^s\tau(t)\,dt\right)\*b.
$$
From the $x$- and $y$-components of the preceding identity, we obtain
$x=c_2(s+c_1)$ and $y=c_3(s+c_1)$, for some constants $c_2,c_3$ satisfying $c_2^2+c_3^2=1$. Therefore, $c_3x=c_2y$, and hence the curve lies in a vertical plane through the $z$-axis.

(3) Write $\*r=u_2\*n+u_3\*b$. Using an argument similar to that in (2), we have
\[
-\kappa u_2=1,
\qquad
u_2'=0,
\qquad
u_3'=u_2+\tau.
\]
Thus $u_2=c_1\neq0$ and
\[
\kappa=-\frac1{c_1},
\qquad
u_3=c_1s+c_2+\int^s\tau(t)\,dt.
\]
The $x$- and $y$-components of $\*r=c_1\*n+u_3\*b$ are
$x=-c_1y'$ and $y=c_1x'$. Therefore
\begin{align*}
x(s)&=A\sin\frac{s}{c_1}+B\cos\frac{s}{c_1},\\
y(s)&=A\cos\frac{s}{c_1}-B\sin\frac{s}{c_1},
\end{align*} for some constants $A, B$. The horizontal arc-length condition gives $A^2+B^2=c_1^2$, and consequently
$x^2+y^2=c_1^2$. Moreover, $xx'+yy'=0$, so the $z$-component is
\[
z(s)=c_1s+c_2+\int^s\tau(t)\,dt.
\]
Hence $\*r$ lies on the circular cylinder $x^2+y^2=c_1^2$, and its $xy$-projection lies on a circle traversed with constant angular speed. Since $z'=c_1+\tau$, the height is affine in $s$ if and only if $\tau$ is constant. If $\tau$ is constant and $c_1+\tau=0$, the curve is a planar circle. If $\tau$ is constant and $c_1+\tau\neq0$, the curve is a nondegenerate circular helix.
\end{proof}

\bibliography{mybib}

@article{blum1966remarkable,
  author = {Blum, Richard},
  title = {A remarkable class of {Mannheim}-curves},
  journal = {Canadian Mathematical Bulletin},
  volume = {9},
  number = {2},
  pages = {223--228},
  year = {1966}
}

@book{cesaro1896lezioni,
  author = {Ces{\`a}ro, Ernesto},
  title = {{Lezioni} di geometria intrinseca},
  publisher = {Presso l'autore-editore},
  year = {1896},
  note = {Reprint by Forgotten Books, 2018}
}

@article{chen2003does,
  author = {Chen, Bang-Yen},
  title = {When does the position vector of a space curve always lie in its rectifying plane?},
  journal = {The American Mathematical Monthly},
  volume = {110},
  number = {2},
  pages = {147--152},
  year = {2003}
}

@article{chiu2018differential,
  author = {Chiu, Hung-Lin and Feng, XiuHong and Huang, Yen-Chang},
  title = {The differential geometry of curves in the {Heisenberg} groups},
  journal = {Differential Geometry and its Applications},
  volume = {56},
  pages = {161--172},
  year = {2018}
}

@article{CHL,
  author = {Chiu, Hung-Lin and Huang, Yen-Chang and Lai, Sin-Hua},
  title = {An application of the moving frame method to integral geometry in the {Heisenberg} group},
  journal = {SIGMA. Symmetry, Integrability and Geometry: Methods and Applications},
  volume = {13},
  pages = {Paper No. 097, 27 pp.},
  year = {2017}
}

@article{choi2013mannheim,
  author = {Choi, Jin Ho and Kang, Tae Ho and Kim, Young Ho},
  title = {{Mannheim} curves in 3-dimensional space forms},
  journal = {Bulletin of the Korean Mathematical Society},
  volume = {50},
  number = {4},
  pages = {1099--1108},
  year = {2013}
}

@book{do2016differential,
  author = {do Carmo, Manfredo P.},
  title = {Differential Geometry of Curves and Surfaces},
  edition = {Revised and updated second edition},
  publisher = {Dover Publications},
  year = {2016}
}

@article{ersoy2012generalized,
  author = {Ersoy, Soley and Tosun, Murat and Matsuda, Hiroo},
  title = {Generalized {Mannheim} curves in {Minkowski} space-time ${E}_1^4$},
  journal = {Hokkaido Mathematical Journal},
  volume = {41},
  number = {3},
  pages = {441--461},
  year = {2012}
}

@article{honda2020bertrand,
  author = {Honda, Shun'ichi and Takahashi, Masatomo},
  title = {{Bertrand} and {Mannheim} curves of framed curves in the 3-dimensional {Euclidean} space},
  journal = {Turkish Journal of Mathematics},
  volume = {44},
  number = {3},
  pages = {883--899},
  year = {2020}
}

@article{jin2008null,
  author = {Jin, Dae-Ho},
  title = {Null {Bertrand} curves in a {Lorentz} manifold},
  journal = {The Pure and Applied Mathematics},
  volume = {15},
  number = {3},
  pages = {209--215},
  year = {2008}
}

@article{lucas2012bertrand,
  author = {Lucas, Pascual and Ortega-Yag{\"u}es, Jos{\'e} Antonio},
  title = {{Bertrand} curves in the three-dimensional sphere},
  journal = {Journal of Geometry and Physics},
  volume = {62},
  number = {9},
  pages = {1903--1914},
  year = {2012}
}

@book{millman1977elements,
  author = {Millman, Richard S. and Parker, George D.},
  title = {Elements of Differential Geometry},
  publisher = {Prentice-Hall},
  address = {Englewood Cliffs, NJ},
  year = {1977}
}

@article{cho2008slant,
  author = {Cho, Jong Taek and Lee, Ji-Eun},
  title = {Slant curves in contact pseudo-{Hermitian} 3-manifolds},
  journal = {Bulletin of the Australian Mathematical Society},
  volume = {78},
  number = {3},
  pages = {383--396},
  year = {2008},
  doi = {10.1017/S0004972708000737}
}

@article{lee2010legendre,
  author = {Lee, Ji-Eun},
  title = {On {Legendre} curves in contact pseudo-{Hermitian} 3-manifolds},
  journal = {Bulletin of the Australian Mathematical Society},
  volume = {81},
  number = {1},
  pages = {156--164},
  year = {2010},
  doi = {10.1017/S0004972709000872}
}

@article{ozgur2012slant,
  author = {{\"O}zg{\"u}r, Cihan and G{\"u}ven{\c{c}}, {\c{S}}aban},
  title = {On some types of slant curves in contact pseudo-{Hermitian} 3-manifolds},
  journal = {Annales Polonici Mathematici},
  volume = {104},
  number = {3},
  pages = {217--228},
  year = {2012},
  doi = {10.4064/ap104-3-1}
}

@article{chiu2019global,
  author = {Chiu, Hung-Lin and Ho, Pak Tung},
  title = {Global differential geometry of curves in three-dimensional {Heisenberg} group and {CR} sphere},
  journal = {The Journal of Geometric Analysis},
  volume = {29},
  number = {4},
  pages = {3438--3469},
  year = {2019},
  doi = {10.1007/s12220-018-00122-x}
}

@article{korpinar2011bertrand,
  author = {K{\"o}rp{\i}nar, Talat and Turhan, Essin},
  title = {Characterization {Bertrand} Curve in the {Heisenberg} Group {Heis}$^3$},
  journal = {International Journal of Open Problems in Complex Analysis},
  volume = {3},
  number = {2},
  pages = {61--67},
  year = {2011}
}

@article{ceylan2015bertrand,
  author = {Ceylan, Ay{\c{s}}e Y{\i}lmaz and Ergin, Abdullah Aziz},
  title = {{Bertrand} mate of a biharmonic curve in {Cartan--Vranceanu} 3-dimensional space},
  journal = {International Electronic Journal of Geometry},
  volume = {8},
  number = {1},
  pages = {45--52},
  year = {2015},
  doi = {10.36890/iejg.592796}
}

@article{tul2011mannheim,
  author = {Tul, S{\i}d{\i}ka and Sar{\i}o{\u{g}}lugil, Ayhan},
  title = {{Mannheim} partner curves in 3-dimensional {Heisenberg} group},
  journal = {International Journal of Physical Sciences},
  volume = {6},
  number = {36},
  pages = {8093--8096},
  year = {2011},
  doi = {10.5897/IJPS10.259}
}

@article{dong2021schwarz,
  author  = {Dong, Yuxin and Ren, Yibin and Yu, Weike},
  title   = {{Schwarz} Type Lemmas for pseudo-{Hermitian} Manifolds},
  journal = {The Journal of Geometric Analysis},
  volume  = {31},
  number  = {3},
  pages   = {3161--3195},
  year    = {2021},
  doi     = {10.1007/s12220-020-00389-z}
}

@article{chang2021sharp,
  title={On the sharp dimension estimate of {CR} holomorphic functions in {Sasakian} manifolds},
  author={Chang, Shu-Cheng and Han, Yingbo and Lin, Chien},
  journal={International Mathematics Research Notices},
  volume={2021},
  number={17},
  pages={12888--12924},
  year={2021},
  publisher={Oxford University Press}
}

@article{balgetir2004null,
  title={Null {Bertrand} curves in {Minkowski} 3-space and their characterizations},
  author={Balgetir, Handan and Bekta{\c{s}}, Mehmet and Inoguchi, Jun-ichi},
  journal={Note di Matematica},
  volume={23},
  number={1},
  pages={7--13},
  year={2004}
}

@article{Lagally1927,
  author  = {Lagally, Max},
  title = {Die {Verwendung} des begleitenden {Dreibeins} f{\"u}r den {Aufbau} der nat{\"u}rlichen {Geometrie}},
  journal = {Sitzungsberichte der Bayerischen Akademie der Wissenschaften, Mathematisch-Naturwissenschaftliche Abteilung},
  year    = {1927},
  pages   = {5--16}
}

@book{Struik,
  author    = {Struik, Dirk J.},
  title     = {Lectures on Classical Differential Geometry},
  edition   = {2nd},
  publisher = {Dover Publications},
  address   = {New York},
  year      = {1988},
  isbn      = {978-0-486-65609-0}
}

@article{okuyucu2017bertrand,
  author  = {Okuyucu, O. Zeki and G{\"o}k, {\.I}smail and Yayl{\i}, Yusuf
             and Ekmekci, Nejat},
  title   = {Bertrand curves in three dimensional {Lie} groups},
  journal = {Miskolc Mathematical Notes},
  volume  = {17},
  number  = {2},
  pages   = {999--1010},
  year    = {2017},
  doi     = {10.18514/MMN.2017.1314}
}

@book{Kowalewski,
  author    = {Kowalewski, Gerhard},
  title = {Vorlesungen {\"u}ber allgemeine nat{\"u}rliche
         {Geometrie} und {Liesche} {Transformationsgruppen}},
  series    = {G{\"o}schens Lehrb{\"u}cherei, I. Gruppe: Reine und angewandte Mathematik},
  volume    = {19},
  publisher = {Walter de Gruyter \& Co.},
  address   = {Berlin and Leipzig},
  year      = {1931}
}

@article{chiu2025invariants,
  author  = {Chiu, Hung-Lin and Lai, Sin-Hua and Liu, Hsiao-Fan},
  title   = {On Invariants of Constant {$p$}-Mean Curvature Surfaces in the {Heisenberg} Group {$H_1$}},
  journal = {SIGMA. Symmetry, Integrability and Geometry:
             Methods and Applications},
  volume  = {21},
  pages   = {Paper No. 011, 25 pp.},
  year    = {2025},
  doi     = {10.3842/SIGMA.2025.011}
}
\end{document}